\documentclass{amsart}
\usepackage[latin1]{inputenc}

\usepackage{subfigure}
\usepackage{mathrsfs}
\usepackage{amscd}
\usepackage{amscd, amssymb, amsmath, amsthm, graphics}
\usepackage{amsmath,amsfonts,amsthm,amssymb}
\usepackage{latexsym,amsmath}
\usepackage{graphicx,psfrag}
\usepackage[all]{xy}
\usepackage{latexsym}
\usepackage{enumerate}
\usepackage{verbatim}
\usepackage{pb-diagram}
\usepackage{color}

\newtheorem{theorem}{Theorem}[section]
\newtheorem{lemma}[theorem]{Lemma}
\newtheorem{prop}[theorem]{Proposition}
\newtheorem{cor}[theorem]{Corollary}

\theoremstyle{definition}
\newtheorem{definition}[theorem]{Definition}

\theoremstyle{remark}
\newtheorem{remark}[theorem]{Remark}

\newtheorem*{claim}{Claim}

\numberwithin{equation}{section}

\newcommand{\CC}{\mathbb{C}}
\newcommand{\RR}{\mathbb{R}}

\newcommand{\ZZ}{\mathbb{Z}}

\newcommand{\Aa}{\mathcal{A}}
\newcommand{\Bb}{\mathcal{B}}

\newcommand{\fF}{\mathcal{F}}
\newcommand{\GG}{\mathcal{G}}

\newcommand{\lL}{\mathscr{L}}

\newcommand{\OO}{\mathscr{O}}

\newcommand{\thmone}{\begin{theorem}}
\newcommand{\thmtwo}{\end{theorem}}
\newcommand{\lemmaone}{\begin{lemma}}
\newcommand{\lemmatwo}{\end{lemma}}
\newcommand{\pfone}{\begin{proof}}
\newcommand{\pftwo}{\end{proof}}
\newcommand{\defone}{\begin{definition}}
\newcommand{\deftwo}{\end{definition}}
\newcommand{\corone}{\begin{cor}}
\newcommand{\cortwo}{\end{cor}}
\newcommand{\cone}{\begin{claim}}
\newcommand{\ctwo}{\end{claim}}
\newcommand{\propone}{\begin{prop}}
\newcommand{\proptwo}{\end{prop}}
\newcommand{\eqone}{\begin{equation}}
\newcommand{\eqtwo}{\end{equation}}
\newcommand{\rmkone}{\begin{remark}}
\newcommand{\rmktwo}{\end{remark}}
\newcommand{\enone}{\begin{enumerate}}
\newcommand{\entwo}{\end{enumerate}}
\newcommand{\itone}{\begin{itemize}}
\newcommand{\ittwo}{\end{itemize}}

\newcommand{\OOO}{\mathscr{O}}

\newcommand{\onehalf}{\left(\begin{array}{cc}}
\newcommand{\theother}{\end{array}\right)}

\DeclareMathOperator{\Diff}{Diff}
\DeclareMathOperator{\Symp}{Symp}
\DeclareMathOperator{\Stab}{Stab}
\DeclareMathOperator{\Fix}{Fix}
\DeclareMathOperator{\id}{id}
\DeclareMathOperator{\Aut}{Aut}
\DeclareMathOperator{\Iso}{Iso}
\DeclareMathOperator{\Map}{Map}
\DeclareMathOperator{\Emb}{Emb}
\newcommand{\CIso}{\mathcal{C}\Iso}
\newcommand{\0}{\mathbf{0}}
\newcommand{\om}{\omega}
\newcommand{\F}{\mathbb{F}}
\newcommand{\Proj}{\mathbb{P}}
\newcommand{\into}{\hookrightarrow}

\newcommand{\oneeq}{\begin{equation}}
\newcommand{\twoeq}{\end{equation}}
\newcommand{\nono}{\noindent}

\newcommand*{\longhookrightarrow}{\ensuremath{\lhook\joinrel\relbar\joinrel\rightarrow}}
\newcommand*{\lthra}{\ensuremath{\relbar\joinrel\twoheadrightarrow}}

\definecolor{vert}{rgb}{0,0.5,0}

\begin{document}

\title[Symplectomorphisms, orbifolds and Lagrangians]{Symplectormophism groups of non-compact manifolds, orbifold balls, and a space of Lagrangians}
\author{Richard Hind, Martin Pinsonnault, Weiwei Wu}

%%%%%%%%%%%%%%%%%%%%%%%%%%%%%%%%%%%%%%%%%%%%%%%%%%%%%%%%%%%%%%%%%%%%%%%%%%%%%%%%
\begin{abstract} We establish connections between contact isometry groups of certain contact manifolds and compactly supported symplectomorphism groups of their symplectizations. We apply these results to investigate the space of  symplectic embeddings of balls with a single conical singularity at the origin.  %and the $\ZZ_n$-equivariant symplectomorphism group of a ball.
Using similar ideas, % by studying automorphisms of symplectizations,
we also prove the longstanding expected result that the space of Lagrangian $\RR P^2$ in $T^*\RR P^2$ is weakly contractible.
\end{abstract}
%%%%%%%%%%%%%%%%%%%%%%%%%%%%%%%%%%%%%%%%%%%%%%%%%%%%%%%%%%%%%%%%%%%%%%%%%%%%%%%%

\maketitle

\noindent {\bf MSC classes:}  53Dxx, 53D35, 53D12


%%%%%%%%%%%%%%%%%%%%%%%%%%%%%%%%%%%%%%%%%%%%%%%%%%%%%%%%%%%%%%%%%%%%%%%%%%%%%%%%
\section{Introduction}
%%%%%%%%%%%%%%%%%%%%%%%%%%%%%%%%%%%%%%%%%%%%%%%%%%%%%%%%%%%%%%%%%%%%%%%%%%%%%%%%

Since the seminal work of Gromov, \cite{Gromov}, the
symplectomorphism groups of closed $4$-manifolds  have been a
subject of much research, see for example \cite{AM}, \cite{Seidel
notes}, as have symplectomorphism groups of manifolds with convex
ends, see for example \cite{Seidel}, \cite{EvansS}. Here we
investigate the simplest symplectic manifolds with both convex and
concave ends, namely the symplectizations $sM$ of $3$-dimensional
contact manifolds $M$. In the case when the contact manifold is a
Lens space $L(n,1)$ the compactly supported symplectomorphism group
$\Symp_c(sL(n,1))$ has a rich topology. In particular, we obtain the
following result:

\thmone\label{theorem:non-compact} The group $\Symp_c(sL(n,1))$,
endowed with  the $C^{\infty}$-topology, has countably many
components, each being weakly homotopy equivalent to the based loop
space of $SU(2)$.  There is a natural map from $\lL(\CIso_n)$, the
based loop group of contact isometry group of $L(n,1)$, to
$\Symp_c(sL(n,1))$ which induces the  weak homotopy equivalence.
\thmtwo

Now, if one of our contact manifolds can be embedded in a
$4$-dimensional  symplectic manifold as a hypersurface of contact
type then there are natural maps from compact subsets of
$\Symp_c(sL(n,1))$ to the symplectomorphism groups of the
$4$-manifold. But as the symplectomorphism group of the $4$-manifold
may have much simpler topology, the induced maps on homotopy groups
will typically be far from injective. For example, $S^3
\hookrightarrow B^4$ as a contact type hypersurface, but while
$\Symp_c(sS^3)$ is weakly homotopy equivalent to the based loop
space of $U(2)$, it is a result of~\cite{Gromov} that $\Symp_c(B^4)$
is contractible.

Our proof of Theorem~\ref{theorem:non-compact} identifies
$\Symp_c(sL(n,1))$ with the based loop space of the K\"ahler
isometry group $K_{n}$ of the Hirzebruch surface
$\F_{n}=\Proj(\OO(n)\oplus\CC)$. Removing the section at infinity
$s_{\infty}$ from $\F_{n}$, and blowing down the zero section
$s_{0}$, one obtains a singular $4$-ball with a conical singularity
of order $n$ at the origin. Since the group $\Symp_c(sL(n,1))$ is
homotopy equivalent to $\Symp_c(\F_{n}\setminus\{s_{\infty}\cup
s_{0}\})$, we can rephrase Theorem~\ref{theorem:non-compact} as a
result on the space of symplectic embeddings of a singular ball of
size $\epsilon\in(0,1)$ into a singular ball of size $1$. In the
second part of the paper, we show that
Theorem~\ref{theorem:non-compact} is equivalent to the following
result:

\thmone\label{thm:orbiballs} The space of symplectic embeddings of a
singular ball of size $\epsilon\in(0,1)$ into a singular ball of
size $1$ is homotopy equivalent to the K\"ahler isometry group
$K_{n}$ of the Hirzebruch surface $\F_{n}$. Moreover, the group of reduced,
compactly supported symplectomorphisms of a singular ball of size $1$ is contractible.
%This implies the compactly supported $\ZZ_n$-equivariant symplectomorphism group of $B(1)$ is contractible.
\thmtwo

Note that in the case $n=1$, the balls are in fact smooth and
Theorem~\ref{thm:orbiballs} reduces to the fact that the space of
symplectic embeddings $B(\epsilon)\into B(1)$ deformation retracts
onto $U(2)$.

In the third part of the paper we apply the techniques used in the proof of Theorem~\ref{theorem:non-compact} in the special case $n=4$ to obtain the homotopy type of a space of Lagrangian submanifolds:

\thmone\label{theorem:lagrangian}
The space of Lagrangian $\RR P^2$ in the cotangent bunble $T^*\RR P^2$, endowed with the $C^{\infty}$-topology, is weakly contractible.
\thmtwo

It is already known that the space of Lagrangian $S^2$ in $T^* S^2$
is contractible, see \cite{Hind}, \cite{Hind2}, and Theorem
\ref{theorem:lagrangian} may be considered as a $\ZZ_2$-equivariant
version.\\

\noindent{\bf Acknowledgements:}  The authors would like to thank the MSRI where part of this work was completed. The second author is supported by a NSERC Discovery Grant. The third author is supported by
NSF Focused Research Grants DMS-0244663.

%%%%%%%%%%%%%%%%%%%%%%%%%%%%%%%%%%%%%%%%%%%%%%%%%%%%%%%%%%%%%%%%%%%%%%%%%%%%%%%%
\section{Symplectomorphism Groups of $sL(n,1)$}\label{section:symplectization}
%%%%%%%%%%%%%%%%%%%%%%%%%%%%%%%%%%%%%%%%%%%%%%%%%%%%%%%%%%%%%%%%%%%%%%%%%%%%%%%%

Consider the lens space
\eqone \label{znaction}
L(n,1) =
\begin{cases}
S^{3} & \text{~} n=1\\
S^{3}/\ZZ_{n} &\text{~}n\geq 2
\end{cases}
\eqtwo
Here, if we think of $S^3$ as the unit sphere in $\CC^2$ then the $\ZZ_n$ action is generated by the diffeomorphism $(z_1,z_2) \mapsto (e^{2\pi i /n} z_1, e^{2\pi i /n} z_2)$.
As contact quotients of $S^3$ with the standard contact form, the lens spaces inherit natural contact one-forms, denoted as $\lambda_{n}$.  There is a standard way to associate a non-compact symplectic manifold to a contact manifold, called the \textit{symplectization}.  Concretely, we consider $L(n,1)\times \RR$ endowed with the symplectic form $d(e^t\lambda_n)$, where $t$ is the coordinate of the second factor $\RR$.  We denote this symplectic manifold $sL(n,1)$. Compactly supported symplectomorphism groups will be denoted by $\Symp_c$.  In this section, we discuss the homotopy type of $\Symp_c(sL(n,1))$, the group of compactly supported symplectomorphisms of $sL(n,1)$.

\subsection{Reducing $sL(n,1)$ to compact manifolds}\label{subsection:reducing to compact manifolds}

We first reduce the problem to the symplectomorphism groups of partially compactified symplectic manifolds.  Let $\OOO(n)$ be the complex
line bundle over $\CC P^1$ with Chern class $c_1=n$.  One can endow the total space of this line bundle with a standard K\"ahler structure, whose
restriction to the zero section is the spherical area form with total area $1$.  We denote the zero section as $C_n$.

\propone\label{prop:liouville and limit} The topological group $\Symp_c(sL(n,1))$ is weakly homotopy equivalent to $\Symp_c(\mathscr{O}(n)\backslash C_n)$.

\pfone Identifying $L(n,1)$ as a circle bundle in $\OO(n)$ with contact structure given by the connection $1$-form, we get a canonical embedding:
\eqone\label{eq:O(n) embedded to symplectization}
\mathscr{O}(n)\backslash C_n\hookrightarrow sL(n,1),
\eqtwo
where the image is $\{(x,t)\in sL(n,1):t<1\}$. Let $\Symp_c^r(sL(n,1))$ be the subgroup of $\Symp_c(sL(n,1))$ consisting of symplectomorphisms supported in $\{t<r\}$, then the embedding (\ref{eq:O(n) embedded to symplectization}) induces an embedding of the corresponding groups of symplectomorphisms, where the image is exactly $\Symp_c^1(sL(n,1))$. On the one hand, for $r>1$, using the inverse Liouville flow one sees that $\Symp_c^r(sL(n,1))$ deformation retracts to $\Symp_c^1(sL(n,1))$; on the other hand, $\Symp_c(sL(n,1))$ is nothing but the direct limit of $\Symp_c^r(sL(n,1))$ as $r\rightarrow\infty$.  This concludes the proof.\pftwo

\proptwo

\subsection{$\Symp_c(\OO(n)\backslash C_n)$ as a loop space}\label{section:Symp as a loop space}

We will find the weak homotopy type of $\Symp_c(\OO(n)\backslash C_n)$
in this section by showing it is weakly homotopy equivalent to a
certain loop space.  We start with some known results about
$\Symp_c(\OO(n))$.

\lemmaone\label{$S_c(O(n))$ is weakly contractible}
$\Symp_c(\OO(n))$ is weakly contractible.
\lemmatwo

This result was shown in \cite{Coffey}, Proposition 3.2. Coffey proceeded by compactifying $\OO(n)$ by adding an infinity divisor to obtain the projectivization of $\OO(n)$, which is the Hirzebruch surface $\F_{n}$. Symplectomorphisms of Hirzebruch surfaces are then studied using holomorphic curves. We note that this can also be deduced from Abreu and McDuff's results in~\cite{AM}.

Now, Coffey also showed that $\Symp_c(\OO(n))$ acts transitively on the space $S(C_{n})$ of unparametrized embedded symplectic spheres in $\OO(n)$ which are homotopic to the zero section. We then have an action fibration
\[
\Stab_c(C_{n})\to \Symp_{c}(\OO(n))\to S(C_{n})
\]
where $\Stab_{c}(C_{n})$ is the subgroup of $\Symp_c(\OO(n))$ consisting of symplectomorphisms which preserve the zero section $C_n$.
\lemmaone[Coffey \cite{Coffey}]\label{lemma:S_c preserving is contractible}
The stabilizer $\Stab_c(C_{n})$, is contractible.
\lemmatwo
Let $\GG_\om(\nu)$ be the symplectic gauge transformations of the normal bundle $\nu$ of $C_{n}$, that is, sections of $Sp(\nu)\to C_{n}$, where $Sp(\nu)$ are the fiberwise symplectic linear maps. Notice that $\GG_{\om}(\nu)\simeq \Map(C_{n},Sp(2))\simeq S^1$ (see~\cite{EvansS}, \cite{Seidel}).

Let $\Fix_c(C_{n})$ be the subgroup of $\Stab_c(C_{n})$ consisting of symplectomorphisms which fix the zero section $C_n$ pointwise. We will use the following lemma from time to time.
\begin{lemma} \label{transitive} The homomorphism $\Fix_c(C_{n}) \to \GG_\om(\nu)$ given by taking derivatives along $C_{n}$ is surjective.
\end{lemma}
\pfone
Let $g \in \GG_\om(\nu)$. Then each $g(z)$ for $z \in S^2$ is a symplectic transformation of the normal fiber $\nu_z$ over $z$. Any such linear symplectic map is the time $1$ Hamiltonian flow $\phi_1$ of a unique quadratic form $Q(z)$ on $\nu_z$.

Consider the Hamiltonian function $H(z,v)= \chi(|v|)Q(z)v$ on $\OO(n)$, where $\chi$ is a bump function equal to $1$ near $0$ and $0$ when $|v|\ge 1$. As $dH=0$ along $C_{n}$ the resulting Hamiltonian flow $\psi_t$ lies in $\Fix_c(C_{n})$. We will check that the corresponding gauge action at time $1$ is precisely $g$.

For this, let $Y \in \nu_z \cong T_0 \nu_z \subset T_z \OO(n)$. Then we claim that $d \psi_t(Y)=\phi_t(Y)$, where in the second term $Y$ is considered as a point in $\nu_z$ and $\phi_t$ is the Hamiltonian flow of $Q: \nu_z \to \RR$. The vector $Y$ can be extended to a Hamiltonian vector field on $\OO(n)$ generated by a function $L$ which is linear on $\nu_z$. Let $X_H$ be the Hamiltonian vector field generated by $H$. Then
\[
{\mathcal L}_{X_H}Y=[X_H,Y]=X_{\{H,L\}}=X_{dH(Y)}
\]
using the same notation throughout for Hamiltonian vector fields. Evaluating at $z$, our Lie derivative is tangent to the fiber $\nu_z$, and restricting to this fiber the function $dH(Y)=dQ(z)(Y)$ is linear and dual under the symplectic form to $X_Q(Y)$. In other words, ${\mathcal L}_{X_H}Y(z)=X_Q(Y)$, identifying two vectors in $\nu_z$. This is equivalent to our claim and so the proof is complete.
\pftwo

Let $\Fix_c^{\id}(C_{n})$ denote the subgroup of $\Fix_c(C_{n})$ consisting of diffeomorphisms whose derivatives act trivially on the normal bundle $\nu$ of the zero section. A simple application of Moser's argument shows that $\Fix_c^{\id}(C_{n})$ is homotopy equivalent to $\Symp_c(\OO(n)\backslash C_n)$, and we will freely switch between these two groups without explicitly mentioning it below.

 Let us write $\Aut_{\om}(\nu)$ for the group of automorphisms of the normal bundle $\nu$ of the zero section $C_{n}$ which are symplectic linear on the fibers and preserve the symplectic form along the zero section. The group $\Stab_c(C_{n})$ acts on $\Aut_{\om}(\nu)$ via its derivative along the zero section. Clearly $\Stab_c(C_{n})$ acts transitively on $C_{n}$ and so by Lemma \ref{transitive} the action on $\Aut_{\om}(\nu)$ is also transitive. Hence we have the fibration
\eqone
\Fix_c^{\id}(C_{n})\longhookrightarrow \Stab_c(C_{n})\lthra \Aut_{\om}(\nu)
\eqtwo

\nono which by Lemma~\ref{lemma:S_c preserving is contractible} yields a weak homotopy equivalence (cf. Proposition 4.66
\cite{Hatcher})
\[
\Fix_c^{\id}(C_{n})\simeq\lL\Aut_{\om}(\nu)
\]
where $\lL\Aut_{\om}(\nu)$ is the space of based loops of $\Aut_{\om}(\nu)$.  Therefore, the following proposition will imply the first part of Theorem \ref{theorem:non-compact}:

\propone\label{prop:loop space of O(n)_l is loop space of SO(3)}
The group $\Aut_{\om}(\nu)$ is homotopy equivalent to the K\"ahler isometry group $K_{n}$ of the Hirzebruch surface $\F_{n}$. In particular,
\[
\Aut_{\om}(\nu)\simeq K_{n}\simeq U(2)/\ZZ_{n}\simeq
\begin{cases}
SO(3)\times S^{1} & \text{~if $n$ is even, $n\neq 0$}\\
U(2) & \text{~if $n$ is odd}
\end{cases}
\]
so that $\lL\Aut_{\om}(\nu)$ has countably many components, where each component is homotopy equivalent to $\lL SU(2)$, that is, to the identity component of $\lL SO(3)$.
\proptwo

\pfone
First notice that $\Aut_{\om}(\nu)$ acts transitively on the symplectic
re\-pa\-ra\-me\-tri\-za\-tion group of the zero section, or equivalently, the symplectomorphism group of $\CC P^1$. We thus have an action fibration
\eqone\label{eq:O(n)_l fibration}
\GG_{\om}(\nu)\longhookrightarrow \Aut_{\om}(\nu)\lthra \Symp(\CC P^1).
\eqtwo
whose fiber is the subgroup which fixes $\CC P^1$ pointwise and thus is simply the gauge group $\GG_\om(\nu)$.

Recall that the Hirzebruch surface $\F_{n}$ is the projectivation $\Proj(\OO(n)\oplus\CC)$. Under the action of its K\"ahler isometry group $K_{n}\simeq U(2)/\ZZ_{n}$, the complex surface $\F_{n}$ is partitionned into three orbits: the zero section $C_{n}$, the section at infinity $C_{n}^{\infty}$ and their open complement $\F_{n}\setminus \{C_{n}\cup C_{n}^{\infty}\}$, see Appendix B in~\cite{AP}. Since the $K_{n}$ action preserves the ruling $\F_{n}\to \CC P^{1}$, every element in $K_{n}$ acts as an isometry of $\CC P^{1}$ and $K_{n}$ acts faithfully on the normal bundle $\nu$ on $C_{n}$ via derivatives. We thus get a commutative diagram of fibrations
\[
\xymatrix{
&\GG_{\om}(\nu) \ar@{->}[r] &\Aut_{\om}(\nu)\ar@{->>}[r] &\Symp(\CC P^1)\\
&S^{1}\ar@{->}[u]\ar@{->}[r] &K_{n}\ar@{->}[u]\ar@{->>}[r] &SO(3)\ar@{->}[u]}
\]
in which the first and third vertical inclusions are homotopy equivalences. It follows that the middle inclusion is a weak homotopy equivalence. Since all spaces involved are homotopy equivalent to CW-complexes, this weak equivalence is a genuine homotopy equivalence. The second part of the statement now follows from substituting $M=SO(3)$ and $N=K_{n}$ in the following simple lemma:

\lemmaone\label{lemma:general algebraic topology} Let $M$ be a
CW-complex with $\pi_2(M)=0$ and $\pi_1(M)$ at most countable.
Suppose $N$ is an $S^1$-bundle over $M$.  Then $\lL(N)$ has
countably many components and we have a weak homotopy equivalence between identity components
$\lL^0(N)\simeq \lL^0(M)$. \lemmatwo

\begin{proof}[Proof of the lemma] This fact is an elementary consequence of the usual
``path-loop" construction. Fix a base point on $N$ and let $P(N) \simeq *$,
be the corresponding based path space. The fibration map
$\pi:N\rightarrow M$ induces the commutative diagram:

\eqone\label{eq:relation between loop of M and N}\xymatrix{
&\lL(N) \ar^{\tilde{\pi}}@{->}[d]\ar@{^{(}->}[r] &P(N)\ar^\pi@{->}[d]\ar@{->>}[r] &N\ar^\pi@{->}[d]\\
&\lL(M)\ar@{^{(}->}[r] &P(M)\ar@{->>}[r] &M }\eqtwo

\nono By assumption, the projections $\pi_*:\pi_k(N)\rightarrow
\pi_k(M)$, are isomorphisms for $k\geq 2$, and the circle fiber and
its multiples are non-zero in $\pi_1(N)$. From the commutative diagram
of the long exact sequence of homotopy groups induced by (\ref{eq:relation
between loop of M and N}), we deduce that:
\eqone\label{eq:higher homotopy group of loop of M and N}
\xymatrix{\tilde{\pi}_*: \pi_k(\lL(N))\ar[r]^-{\cong}&\pi_k(\lL(M))}
, \hskip 3mm \text{when} \hskip1mm k\geq 1;
\eqtwo

\nono Moreover, we have noticed that $\pi_1(N)$ is the central
extension of $\ZZ$ and $\pi_1(M)$, hence the lemma follows.
%\eqone\label{eq:components of loop of M and N} 0\rightarrow
%\pi_1(S_c(\OO(n))_{\bar{0}})\rightarrow\pi_2(SO(3))\rightarrow
%\pi_1(S^1)\rightarrow \pi_0(S_c(\OO(n))_0)\rightarrow 0. \eqtwo
\end{proof}

\begin{comment}
\eqone\label{eq:Sc(O(n)) homotopy type}
\pi_k(S_c(\OO(n))_{\bar{0}})=\pi_k(S_c(\OO(n))_0)=\pi_{k+1}(SO(3)), \hskip 3mm \text{when} \hskip1mm k\geq 2;
\eqtwo
\nono and,
\eqone\label{eq:Sc(O(n)) remaining part}
0\rightarrow \pi_1(S_c(\OO(n))_{\bar{0}})\rightarrow\pi_2(SO(3))\rightarrow \pi_1(S^1)\rightarrow \pi_0(S_c(\OO(n))_0)\rightarrow 0.
\eqtwo
\end{comment}
\nono This concludes the proof of Proposition~\ref{prop:loop space of O(n)_l is loop space of SO(3)}
\pftwo

\subsection{The loop group of the contact isometries of $L(n,1)$} \label{contisom}

In this section, we prove the second part of Theorem \ref{theorem:non-compact} by showing that a natural inclusion map is a weak homotopy equivalence.
Unlike the usual notion of contactomorphism which preserves only the
contact structures, we need to consider the automorphisms of
$L(n,1)$ called \textit{contact isometries}. These are diffeomorphisms which preserve the contact form $\lambda_n$ and the round metric induced from the round metric on $S^3$ under projection. We denote the group of contact isometries of the lens spaces of $L(n,1)$ as $\CIso_{n}$.  It acts on $L(n,1)$ in such a way that the Reeb orbits are preserved. Therefore, if we think of $L(n,1)$ as a unit circle bundle in $\OO(n)$ with the Reeb orbits as the circle fibers, there is an induced isometric action of $\CIso_{n}$ on the base $\CC P^1$ endowed with the standard round metric. Also, since the action on the fibers is linear, there is a natural inclusion $\CIso_{n} \hookrightarrow \Aut_{\om}(\nu)$. Therefore, along with (\ref{eq:O(n)_l fibration}), one obtains the following diagram of fibrations :
\eqone\label{eq:contact isometry}
\xymatrix{
&S^1 \ar@{^{(}->}[d]\ar@{^{(}->}[r] &\CIso_{n}\ar@{^{(}->}[d]\ar@{->>}[r] &SO(3)\ar@{^{(}->}[d]\\
&\GG_{\om}(\nu)\ar@{^{(}->}[r] &\Aut_{\om}(\nu)\ar@{->>}[r] &\Symp(\CC P^1) }
\eqtwo

\nono Notice that we have weak homotopy equivalences in both the base and fiber.  Therefore, the natural inclusion of $\CIso_{n}$ into $\Aut_{\om}(\nu)$ is in fact a (weak) homotopy equivalence.

We now want to describe a natural map from $\lL(\CIso_{n})$ to $\Symp_c(\OO(n)\backslash C_n)$ (or equivalently $\Fix^{\id}_{c}(C_{n})$, see section~\ref{section:Symp as a loop space}) which induces a weak homotopy equivalence. Given Proposition \ref{prop:liouville and limit} this will imply the remainder of Theorem \ref{theorem:non-compact}. To this end, consider the smooth path space
\[
P(\CIso_{n})=\{\phi:(-\infty,+\infty)\rightarrow \CIso_{n}: \phi(t)=id, t\leq0, \phi(t)=\phi(1), t\geq1\}.
\]

This is just the usual based path space when restricted to $t\in[0,1]$, thus it is a contractible space. Given $\phi\in P(\CIso_{n})$ one can define the following diffeomorphism of $sL(n,1)$:

\begin{align}
\phi':L(n,1)\times \RR&\longrightarrow L(n,1)\times \RR\\
           (x,t)      &\longmapsto     (\phi(t)x,t)
\end{align}

By definition, $\phi'|_{t\leq0}=id$, and $\phi'|_{1\leq t\leq 2}$ is a symplectomorphism
induced by a contact isometry multiplied by identity in the $\RR$-direction.  However $\phi'$ fails to be a symplectomorphism in general. Let $\omega_0 = d(e^t \lambda_n)$, the canonical symplectic form on $sL(n,1)$, and $\omega_1=\phi'^* \omega_0$. Then nevertheless we claim that the exact forms
$\omega_u=(1-u)\omega_0 + u\omega_1$ are symplectic for all $0 \le u \le 1$.

{\it Proof of claim.} To see this, arguing by contradiction, note that if an $\omega_u$ fails to be symplectic then it has a kernel of dimension at least $2$, which must intersect the tangent space to some level $L(n,1) \times \{t\}$ nontrivially. As our $\phi(t)$ are contact isometries this kernel must be the kernel of $d\lambda_n$, namely the Reeb direction. But as the Reeb direction is preserved and $\phi'_*(\frac{\partial}{\partial t})$ always has a positive
$\frac{\partial}{\partial t}$ component, the Reeb vector pairs nontrivially with $\frac{\partial}{\partial t}$ under all $\omega_u$.\\

Given our claim, we can apply Moser's method, see \cite{MSI}, that is, we compose $\phi'$ with the flow of the time-dependent vector field defined by $X_u \lfloor \omega_u = \phi'^*(e^t \lambda_n) - e^t \lambda_n$. The composition is
a symplectomorphism $\tilde{\phi}$ of $sL(n,1)$ supported in $\{t \ge 0\}$. On $\{t \ge 1\}$ we have $\phi'^*(e^t \lambda_n) = e^t \lambda_n$ and so on this region we have $\tilde{\phi}|_{t\geq1} = \phi'|_{t\geq1}$.

Next, as $\phi'$ is translation invariant on $\{t \ge 1\}$ we can perform a symplectic cut at the level of $\{t=2\}$. Recall from Section \ref{subsection:reducing to compact manifolds} that, up to a scaling, we may identify symplectically $\OO(n)\backslash C_n$ with $sL(n,1)_{t<2}$. Hence $\tilde{\phi}$ descends to a compactly supported symplectomorphism of $\OO(n)$ preserving the zero-section, that is, we have a map $P(\CIso_{n}) \to \Stab_c(C_{n})$, $\phi \mapsto \tilde{\phi}$.\\

%Observe that, given a cylindrical almost complex structure $J$ which is tamed by $\omega_n=d(d^t\lambda_n)$,
%$\phi'^*\omega_n$ is still tamed by the same $J$.  Therefore, Moser's trick yields a symplectomorphism $\tilde{\phi}\in Symp_c^2(sL(n,1))$.
%Since $\tilde{\phi}|_{t\geq1}=\phi'|_{t\geq1}$, where Moser's flow vanishes, one is allowed to

\nono{\bf Claim:} The following diagram of fibrations is commutative and all maps are continuous.  The rightmost vertical arrow is a homotopy equivalence:

\eqone\label{eq:path-loop of Cont}\xymatrix{
&\lL(\CIso_{n})\ar^{}@{->}[d]\ar@{^{(}->}[r] &P(\CIso_{n})\ar^{}@{->}[d]\ar@{->>}[r] &\CIso_{n}\ar@{_{(}->}[d]\\
&\Fix_c^{\id}(C_{n})\ar@{^{(}->}[r] &\Stab_c(C_{n})\ar@{->>}[r] &\Aut_{\om}(\nu)}
\eqtwo

\pfone[Proof of claim:] The second arrow of the first row is simply the restriction of an element $\phi$ to $\phi(2)$.  The continuity of the vertical maps follows from the continuous dependence of solutions of an ODE on initial conditions when applying Moser's method.  The rightmost vertical arrow is the one induced from~(\ref{eq:contact isometry}) and is thus a homotopy equivalence.  The commutativity of the diagram (\ref{eq:path-loop of Cont}) is straightforward  from definitions.
\pftwo

Now the middle vertical arrow is a homotopy equivalence due to the contractibility of both spaces, see Lemma \ref{lemma:S_c preserving is contractible}, and the rightmost arrow is also a weak homotopy equivalence from the argument at the start of this subsection. Therefore, the leftmost vertical arrow is a homotopy equivalence as well, and provides the desired mapping.  Hence, the second part of Theorem \ref{theorem:non-compact} follows.

\begin{comment}
\nono Therefore, the long exact sequence gives

\[
\pi_k(\CIso_{n})=\pi_k(SO(3)), \hskip 3mm k\geq 2,
\]
\[
0\longrightarrow \ZZ\longrightarrow \pi_1(\CIso_{n})\longrightarrow \ZZ_2\longrightarrow 0.
\]

In either case of extension for second short exact sequence, we see that the loop space $\lL(\CIso_{n})$ is weak homotopy equivalent to $\Symp_c(\OO(n)\backslash C_n)$ in components, and both of them have $\ZZ$ copies of connected components. Therefore, ...

\end{comment}

%%%%%%%%%%%%%%%%%%%%%%%%%%%%%%%%%%%%%%%%%%%%%%%%%%%%%%%%%%%%%%%%%%%%%%%%%%%%%%%%
\section{Space of Symplectic Embeddings of Orbifold Balls}
%%%%%%%%%%%%%%%%%%%%%%%%%%%%%%%%%%%%%%%%%%%%%%%%%%%%%%%%%%%%%%%%%%%%%%%%%%%%%%%%

In this section, we study the space of symplectic embeddings of balls with a single conical singularity at the origin. We first briefly recall the two different notions of maps between orbifolds that we use and the related definition of orbifold embeddings. A comprehensive discussion of orbifold structures and of orbifold maps can be found in~\cite{BB}.

Given an orbifold $\Aa$, we write $|\Aa|$ for its underlying topological space. An \emph{unreduced} orbifold map $(f,\{\hat{f}\})$ between two orbifolds $\Aa$ and $\Bb$ consists of the following data:
\begin{enumerate}
\item a continuous map $f:|\Aa|\to|\Bb|$ of the underlying topological spaces;
\item for all $x\in|\Aa|$, the choice of a germ $\hat{f}_{x}$ of local lift of $f$ to uniformizing charts $U$ and $V$ centered at $x$ and $f(x)$.
\end{enumerate}
A \emph{reduced} orbifold map is a continuous map $f:|\Aa|\to|\Bb|$ of the underlying topological spaces such that smooth lifts exist at every point.
%is a homeomorphism $|f|:|\Aa|\to|\Bb|$ such that, given any point $a\in\Aa$, we can find uniformizing charts $U$ and $V$ centered at $a$ and $|f|(a)$, together with a smooth equivariant lift $\hat{f}_{a}:U\to V$, such that the following diagram commutes
%\[
%+++
%\]
The set of smooth unreduced orbifold maps between $\Aa$ and $\Bb$ will be denoted by $C_{orb}^{\infty}(\Aa,\Bb)$, while we will write $C_{red}^{\infty}(\Aa,\Bb)$ for the set of smooth reduced orbifold maps. Smooth unreduced or reduced diffeomorphisms are defined accordingly by requiring $f$ to be a homeomorphism and all lifts to be smooth local diffeomorphisms. The sets of all unreduced or reduced diffeomorphisms of an orbifold $\Aa$ can be naturally endowed with a $C^{\infty}$ topology that make them Fréchet Lie groups. The short exact sequence
\[
1\to\Gamma_{\id}\to\Diff_{orb}(\Aa)\to\Diff_{red}(\Aa)\to 1
\]
is then a principal bundle whose fiber $\Gamma_{\id}$ is the (discrete) group of all lifts of the identity map.

A smooth (unreduced, resp. reduced) embedding $f:\Aa \to \Bb$ is a smooth (unreduced, resp. reduced) orbifold map which is a diffeomorphism onto its image and which covers a topological embedding $f:|\Aa| \to |\Bb|$. We will denote by $\Emb_{orb}(\Aa,\Bb)$ and $\Emb_{red}(\Aa,\Bb)$ the corresponding embedding spaces.

If all the uniformizing charts are symplectic, and if all the local group actions preserve the symplectic forms, the orbifold atlas is said to be  symplectic. Orbifold symplectic maps are then defined in the obvious way. In particular, since the open set of regular points $\Aa_{reg}$ becomes an open symplectic manifold, orbifold symplectomorphisms restrict to genuine smooth symplectomorphisms of $\Aa_{reg}$.

Let us write $B_{n}(\epsilon)$ for a symplectic ball of size $\epsilon$ with a single conical singularity of order $n\geq 1$ at the origin, that is,
\[
B_{n}(\epsilon) := B^{4}(n\epsilon)/\ZZ_{n}
\]
where $B^{4}(r)$ stands for the standard ball of radius $\sqrt{r/\pi}$ in $\RR^{4}$ and $\ZZ_n$ acts as in (\ref{znaction}). In this paper, we are only interested in the simplest possible embedding spaces between symplectic orbifolds, namely $\Emb_{red}(B_{n}(\epsilon), B_{n}(1))$. In that case, it is easy to see that $\Emb_{orb}(B_{n}(\epsilon), B_{n}(1))$ consists of smooth symplectic embeddings $f:B^{4}(n\epsilon)\to B^{4}(n)$ of standard smooth balls that are equivariant with respect to the standard $\ZZ_{n}$ action, and that
\[
\Emb_{red}(B_{n}(\epsilon), B_{n}(1))=\Emb_{orb}(B_{n}(n\epsilon), B_{n}(n))/\ZZ_{n}
\]
This follows from the fact that any local smooth lift at the conical point extends uniquely to the whole ball, see~\cite{BB}. Since the space of smooth symplectic embeddings retracts onto $U(2)$, and since the $\ZZ_{n}$ action belongs to the center of $U(2)$, one can show that the space of $\ZZ_{n}$-equivariant embeddings of smooth balls is itself homotopy equivalent to $U(2)$, see~\cite{VH}. Therefore, we have the following results:
\begin{prop}\label{reduced_embeddings}
The space of reduced symplectic embeddings $\Emb_{red}(B_{n}(\epsilon), B_{n}(1))$ is homotopy equivalent to
\[K_{n}:=U(2)/\ZZ_{n}\simeq
\begin{cases}
SO(3)\times S^{1} & \text{~if $n$ is even, $n\neq 0$}\\
U(2) & \text{~if $n$ is odd}
\end{cases}
\]
\end{prop}

Just as in the smooth case, one can show that the group of compactly supported and reduced symplectomorphisms of the open orbifold ball $B_{n}(1)$ acts transitively on $\Emb_{red}(B_{n}(\epsilon), B_{n}(1))$, see~\cite{VH}. We get an action fibration
\[
\Stab_{c,red}(\iota) \to \Symp_{c,red}(B_{n}(1)) \to \Emb_{red}(B_{n}(\epsilon), B_{n}(1))
\]
where $\iota:B_{n}(\epsilon)\to B_{n}(1)$ is the inclusion, and where $\Stab_{c,red}(\iota)$ is the subgroup made of those reduced symplectomorphisms that are the identity on the image $\iota(B_{n}(\epsilon))$. This subgroup is homotopy equivalent to the group of reduced diffeomorphisms that are the identity near the image $\iota(B_{n}(\epsilon))$. Performing a symplectic blow-up of the ball $\iota(B_{n}(\epsilon))$, those symplectomorphisms lift to symplectomorphisms of the Hirzebruch surface $\F_{n}$ that are the identity near the zero section and near the section at infinity. This last group is itself homotopy equivalent to $\Symp_c(sL(n,1))$. Hence, we get a homotopy fibration
\[
\Symp_c(sL(n,1))\to \Symp_{c,red}(B_{n}(1)) \to \Emb_{red}(B_{n}(\epsilon), B_{n}(1))
\]
which shows that the homotopy equivalence $\Symp_c(sL(n,1))\simeq \lL K_{n}$, together with Proposition~\ref{reduced_embeddings}, imply the following mild generalization of a fundamental result due to Gromov:
\begin{prop}
The group $\Symp_{c,red}(B_{n}(1))$ of reduced, compactly supported symplectomorphisms of an open ball of size $1$ with a single conical singularity of order $n$ at the origin is contractible.
\end{prop}
This completes the proof of Theorem~\ref{thm:orbiballs}.

%%%%%%%%%%%%%%%%%%%%%%%%%%%%%%%%%%%%%%%%%%%%%%%%%%%%%%%%%%%%%%%%%%%%%%%%%%%%%%%%
\section{Space of Lagrangian $\RR P^2$ in $T^*\RR P^2$}
%%%%%%%%%%%%%%%%%%%%%%%%%%%%%%%%%%%%%%%%%%%%%%%%%%%%%%%%%%%%%%%%%%%%%%%%%%%%%%%%

We prove Theorem \ref{theorem:lagrangian}.  Let the space of Lagrangian $\RR P^2$ in $T^*\RR P^2$ be denoted as $\mathcal{L}$. The group of compactly supported Hamiltonian symplectomorphisms of $T^{*}\RR P^{2}$ acts transitively on $\mathcal{L}$, see~\cite{Hind2}, and our point of departure is the corresponding action fibration
\eqone\label{eq:main fibration}
\Stab_{c}(\0)\longhookrightarrow \Symp_c(T^*\RR P^2)\lthra\mathcal{L}
\eqtwo
where $\Stab_{c}(\0)$ is the subgroup of $\Symp_c(T^*\RR P^2)$ which preserves the zero section~$\0$.

Notice that for any Lie group $G$, $\pi_0(G)$ inherits a natural group structure from $G$. It is proved in \cite{EvansS} that:

\thmone
 $\Symp_c(T^*\RR P^2)$ is weakly homotopic to $\ZZ$. Moreover, the generator of $\pi_0(\Symp_c(T^*\RR P^2))$ as a group is the generalized Dehn twist in $T^*\RR P^2$. \thmtwo

We will also make use of the following fact, which may be well-known but for which the authors unfortunately know of no reference:

\lemmaone \label{group} Let $H\longhookrightarrow G\lthra B$ be a homotopy fibration where $H\lhd G$ are groups.  Then the following two maps in the induced long exact sequence are both group homomorphisms:
\[
\xymatrix{\pi_1(B)\ar[r]^i &\pi_0(H)\ar[r]^j &\pi_0(G) }
\]
\lemmatwo

\pfone Let $x_0\in B$ be the image of $id\in G$.  Given a loop $\alpha:[0,1]\rightarrow B$, $\alpha(0)=\alpha(1)=x_0$, let $\bar{\alpha}$ be the lift of $\alpha$ and $i(\alpha)$ be the connected component of $H$ where $\bar{\alpha}(1)$ lies. Consider another loop $\beta:[0,1]\rightarrow B$, $\beta(0)=\beta(1)=x_0$, then the lift of concatenation $\alpha\#\beta$ can be chosen to be

\[ \label{eq:1} \overline{\alpha\#\beta}(t)=\left\{ \begin{aligned}
         \bar{\alpha}(2t), \hskip 2mm t\leq \frac{1}{2} \\
                  \bar{\alpha}(1)\cdot\bar{\beta}(2t-1), \hskip 2mm
                  t>\frac{1}{2}
                          \end{aligned} \right.
\]

Therefore,
$i(\alpha\#\beta)=\overline{\alpha\#\beta}(1)=\bar{\alpha}(1)\bar{\beta}(1)$,
verifying the claim for the map $i$.  The fact that $j$ is a homomorphism is trivial because the inclusion $H\hookrightarrow G$ is a homomorphism. \pftwo

To compute the homotopy type of $\Stab_{c}(\0)$ we need to consider the diffeomorphism group of $\RR P^2$.  We have the following result, of which the proof is postponed to the appendix:

\propone\label{prop:Diffrp2} The diffeomorphism group of $\RR P^2$ is weakly homotopic to $SO(3)$. Moreover, the standard inclusion is a weak homotopy equivalence.\proptwo

With this understood, we define $\Fix_{c}(\0)$ to be the subgroup of compactly supported symplectomorphisms of $T^*\RR P^2$ which fixes the zero section pointwise. We obtain a further action fibration:
\[
\Fix_{c}(\0)\longhookrightarrow \Stab_{c}(\0)\lthra \Diff(\RR P^2).
\]
We may also consider the following object: given the standard round metric $g_0$ on $\RR P^2$, let $\Stab_{c}^{\Iso}(\0)$ be the symplectomorphisms which are compactly supported and induce an isometry on the zero section.

\begin{comment}
The homotopy fiber, as is shown in the Section~\ref{section:symplectization}, is weakly homotopic to countably many copies of $SO(3)$, moreover, the $\pi_0$ inherits a group structure of $\ZZ\oplus\ZZ/2\ZZ$ or $\ZZ$ from the group structure of the symplectomorphism groups. Therefore, the above action-fibration gives, for $k\geq1$,
%
\[
\pi_{k+1}(\Fix_{c}(\0))\rightarrow \pi_{k+1}(\Diff(\RR P^2))\overset{\cong}{\rightarrow}\pi_k(S_{\overline{\RR P^2}}).
\]
%
This implies $\pi_{k}(\Fix_{c}(\0))=0$ for $k\geq 2$.  For the last couple terms of the long exact sequence, we have:
%
\[
0\rightarrow\pi_{1}(\Fix_{c}(\0))\rightarrow \pi_{1}(\Diff(\RR P^2))\rightarrow\pi_0(S_{\overline{\RR P^2}})\rightarrow\pi_0(\Fix_{c}(\0))\rightarrow0.
\]
%
The reasonable guess is that the connecting map $\pi_1(\Diff(\RR P^2))\rightarrow \pi_0(S_{\overline{\RR P^2}})$ is injective, so that
this leaves $\Fix_{c}(\0)$ homotopic equivalent to countably many points.  Inserting this fact back to (\ref{eq:main fibration}), we obtain:

\end{comment}

Now we have the following commutative diagram of fibrations:
\eqone\xymatrix{
&\Fix_{c}(\0) \ar@{=}[d]\ar@{^{(}->}[r] &\Stab_{c}^{\Iso}(\0)\ar@{_{(}->}[d]\ar@{->>}[r] &SO(3)\ar@{_{(}->}[d]\\
&\Fix_{c}(\0)\ar@{^{(}->}[r] &\Stab_{c}(\0)\ar@{->>}[r] &\Diff(\RR P^2)
}\eqtwo
From Proposition \ref{prop:Diffrp2}, we observe that the vertical
arrows on the two sides are weak homotopy equivalences, so the middle one
is also a weak homotopy equivalence.  Note also that the inverse Liouville flow
contracts $T^*\RR P^2$ to the zero section.  Now, taking into consideration the bundle
metric on $T^*\RR P^2$ induced by $g_0$, we may talk about the length of cotagent vectors.
By the same direct limit and Liouville flow argument as in Proposition \ref{prop:liouville
and limit} we may restrict our attention to the symplectomorphisms
supported in $T^*_r\RR P^2$, which consists of cotagent vectors with length $\leq r$ for some $r>0$.
We will assume that $r=1$ below.

\lemmaone\label{lemma:S_rigid preserving is Z}
\begin{enumerate}[(i)]
\item $\Stab_{c}^{\Iso}(\0)$ is weakly homotopy equivalent to
$\ZZ$;
\item $\pi_0(\Stab_{c}^{\Iso}(\0))$ is isomorphic to $\ZZ$ as a group.
\end{enumerate}
\lemmatwo

\begin{remark}
It is very tempting to conclude $(i)$ directly from the results in the previous sections by setting $n=4$, see the first paragraph of the proof.  However, the connecting map in $(4.2)$ seems then difficult to understand directly. That is why we use a slightly different argument below.
\end{remark}

\begin{proof}[of Lemma~\ref{lemma:S_rigid preserving is Z}] We first notice the following fact: a symplectomorphism which fixes a smooth Lagrangian pointwise also fixes the framing of the Lagrangian.  This follows from the corresponding linear statement that, a symplectomorphism of $T^*M$ which is linear on the fibers is indeed a cotangent map of a diffeomorphism on the base.  This is also used in~\cite{Coffey}, proof of Theorem 1.3. It follows from this that the subgroup of $\Stab_{c}^{\Iso}(\0)$ consisting of maps which act on a neighborhood of $\RR P^2$ by the cotangent map of an isometry of $\RR P^2$ is weakly homotopy equivalent to $\Stab_{c}^{\Iso}(\0)$. Therefore we are able to consider this subgroup instead of $\Stab_{c}^{\Iso}(\0)$, and use the same notation to denote it throughout the rest of the proof.

Given $\psi\in \Stab_{c}^{\Iso}(\0)$, denoting the cotangent map of $\psi|_{\RR P^2}$ as $c_\psi$, we may consider the symplectomorphism $\tilde{\psi}:=c_\psi^{-1}\circ\psi$ on $T^*\RR P^2$. The map $\tilde{\psi}$ is not compactly supported in $T^*\RR P^2$ , but it fixes $\RR P^2$ pointwise and thus (by our assumption that the maps are cotangents near $\RR P^2$) also a neighborhood.  Since $\psi$ preserves the round metric on $\RR P^2$, $c_\psi$ preserves the Reeb vector field on the level sets of $T^*\RR P^2$.  Therefore, by a symplectic cut on the level set $r=1$, one obtains a symplectomorphism $\psi'$ of $T^*\RR P^2$ cut along the level $r=1$. This symplectic manifold is just $\CC P^2$ with the standard symplectic form $\omega_{FS}$, see \cite{Audin} and~\cite{Lerman}.  From the construction, $\psi'$ preserves the symplectic reduction of the boundary, a symplectic $(+4)$-sphere which is indeed the quadratic sphere $\{[x,y,z]\in\CC P^2: x^2+y^2+z^2=0\}$, and it fixes a neighborhood of the standard Lagrangian $\RR P^2=Re(\CC P^2)$.  Removing the Lagrangian $\RR P^2$, one sees that $\psi'$ descends to a compactly supported symplectomorphism of $\OO(4)$, which is denoted as $\bar{\psi}$. Define $H$ to be the image of the bar assignment $\psi\mapsto\bar{\psi}$, $H$ is clearly homeomorphic to $\Stab_{c}^{\Iso}(\0)$. In the rest of the proof we investigate the homotopy type of $H$.

\lemmaone\label{lemma:S acts transitively on zero section}
Let $U$ be a sufficiently small neighborhood of the zero section in $\OO(4)$, then $H|_U=SO(3)$, and it acts transitively on the zero section.
\lemmatwo

\pfone From the construction, there is a surjective map $f:SO(3)\rightarrow H|_U$.  But remembering that points of the zero section in $\OO(4)$ corresponds to the lifts of a geodesic, the map $f$ is clearly injective.
\pftwo

\rmkone\label{rmk:Liu Yi's construction}
There is an interesting model described to the authors by Yi Liu. Consider $\RR^3$ with the standard Euclidean metric $g_E$. Consider an oriented normal frame $(e_1,e_2,e_3)$ as a point on $\RR P^3$, it fibers over $S^2$ by projection to $e_1$.  Let $\varpi:\RR P^3\rightarrow \RR P^3$ be the involution sending $(e_1,e_2,e_3)$ to $(-e_1,-e_2,e_3)$ and consider its quotient $L(1,4)$. This can be identified with the unit cotengent bundle of $\RR P^{2}$ and fibers over $\RR P^2$ by the projection
\eqone\label{eq:L(1,4) projects to
RP^2}[e_1,e_2,e_3]\rightarrow[e_1]
\eqtwo
with fiber $S^1$.  On the other hand, one may project $L(1,4)$
to $S^2$ by sending
\eqone\label{eq:L(1,4) projects to S^2}[e_1,e_2,e_3]\rightarrow e_3
\eqtwo
Endow all the spaces involved with the metric inherited by $g_E$ and
use the obvious $SO(3)$ action on $\RR P^2$ as constructed, then the
projections interplay correctly with the symplectic structure on
$T^*\RR P^2$. In other words, given an isometry of $\RR P^2$,
represented by an element $R\in SO(3)$ in the above model,
the corresponding action on the unit cotangent bundle is
described by the same element $R$ acting on $L(1,4)$. In turn $R$ acts
on the fibration (\ref{eq:L(1,4) projects to S^2}). In this way we
retrieve the action of $H|_U$.

\rmktwo

\nono Returning to the proof of Lemma~\ref{lemma:S_rigid preserving is Z}, given the round metric $g$ on the zero section $C_{4}$ of $\OO(4)$, we
consider the subgroups
\begin{align*}
\Stab_c(C_{4})&=\{\psi\in \Symp_c(\OO(4))~:~\psi\text{~preserves the zero section~} C_{4} \} \\
\Stab_c^{\Iso}(C_{4})&=\{\psi\in \Stab_c(C_{4})~:~\psi \text{ restricted to the zero section}\\
 &\text{\hspace*{2.3cm} is an isometry with respect to the metric }g \}\\
\Fix_{c}(C_{4})&=\{\psi\in\Stab_c^{\Iso}(C_{4})~:~\psi|_{C_{4}} = \id \}
\end{align*}

\nono We then again have a commutative diagram of fibrations:
\eqone\xymatrix{
&\Fix_{c}(C_{4}) \ar@{=}[d]\ar@{^{(}->}[r] &\Stab_c^{\Iso}(C_{4})\ar@{_{(}->}[d]\ar@{->>}[r] &SO(3)\ar@{_{(}->}[d]\\
&\Fix_{c}(C_{4})\ar@{^{(}->}[r] &\Stab_c(C_{4})\ar@{->>}[r] &\Symp(\CC P^1)
}\eqtwo
Using the fact that the embedding of $SO(3)$ into $\Symp(\CC P^1)$ is a weak homotopy equivalence and Lemma \ref{lemma:S_c preserving is contractible}, we deduce that $\Stab_c^{\Iso}(C_4)$ is also weakly contractible.

Now consider the subgroup $H\subset \Stab_c^{\Iso}(C_4)$. We construct the following group homomorphism from $\Stab_c^{\Iso}(C_4)$ to the gauge group:
\[
\phi: \Stab_c^{\Iso}(C_4)\longrightarrow \Map(S^2,Sp(2))\simeq S^1
\]
To define $\phi$, let $t\in \Stab_c^{\Iso}(C_4)$.  Then $t|_{C_4}$ acts on the zero section $C_4$ isometrically. By Lemma \ref{lemma:S acts transitively on zero section} and the remark following it, there exists an element $u\in H$ such that $u|_{C_4}=t|_{C_4}$. Now we define $\phi(t)$ to be the gauge of $t\cdot u^{-1}$.

Notice first that for any $u\in H$, its action on the normal bundle of $C_4$ is uniquely determined by its action on $C_4$, hence $\phi(t)$ does not depend on the choice of $u$ and is well-defined.

The homomorphism $\phi$ is clearly surjective by Lemma \ref{transitive} since $\Fix_{c}(C_4)\subset\Stab_c^{\Iso}(C_4)$.

It is also not hard to verify that $\ker(\phi)\simeq H$:  indeed, for $t\in \ker(\phi)$, by definition there exists $u\in H$, such that $t\cdot u^{-1}$ acts trivially on the normal bundle of $C_{4}$. However, up to homotopy, $\ker(\phi)$ consists of $t$ for which there is a $u\in H$ such that $t\cdot u^{-1}$ acts trivially on a neighborhood of $C_{4}$. As all elements acting trivially on a neighborhood lie in $H$, we deduce that all such $t$ lie in $H$ too.

Therefore we have the following fibration:
%
\begin{comment}
By Lemma \ref{lemma:S acts transitively on zero section} and the remark following it, we see the homogeneous space of $H$ is the full gauge group $Map(S^2,S^1)\sim S^1$, because the action of $H$ on a neighborhood $U$ of the zero section is completely determined by the action on the zero section.
\end{comment}
%
\[
H\longhookrightarrow \Stab_c^{\Iso}(C_4)\lthra S^1
\]
which implies that $H$ is weakly homotopy equivalent to $\ZZ$ and, by Lemma \ref{group}, that $\pi_0(S)\cong \ZZ$ since $\Stab_c^{\Iso}(C_4)$ is
contractible. This concludes our proof of Lemma~\ref{lemma:S_rigid preserving is Z}.
\end{proof} % of lemma:S_rigid preserving is Z

\pfone[Proof of Theorem \ref{theorem:lagrangian}:] For $\pi_i(\lL)$, $i\geq1$ the theorem follows immediately from Lemma \ref{lemma:S_rigid preserving is Z} and the homotopy fibration (\ref{eq:main fibration}).  Since the Dehn twists are also contained in the subgroup $\Stab_{c}^{\Iso}(\0)$, one sees that the map
$\pi_0(\Stab_{c}(\0))\rightarrow \pi_0(\Symp_c(T^*\RR P^2))$ is surjective.  However, by Lemma \ref{group}, since both groups are $\ZZ$, it can only be an
isomorphism.
\pftwo

\newpage
\appendix
%%%%%%%%%%%%%%%%%%%%%%%%%%%%%%%%%%%%%%%%%%%%%%%%%%%%%%%%%%%%%%%%%%%%%%%%%%%%%%%%
\section{The Diffeomorphism Group of $\RR P^2$}
%%%%%%%%%%%%%%%%%%%%%%%%%%%%%%%%%%%%%%%%%%%%%%%%%%%%%%%%%%%%%%%%%%%%%%%%%%%%%%%%

%\centerline{\Large {\bf Appendix} }\medskip
%\centerline{\Large Diffeomorphism group of $\RR P^2$}
%\medskip
\begin{figure}[b]
 \centering
 \subfigure[``Example'' of a curve in $\fF_2$]{
  \includegraphics[scale=0.4]{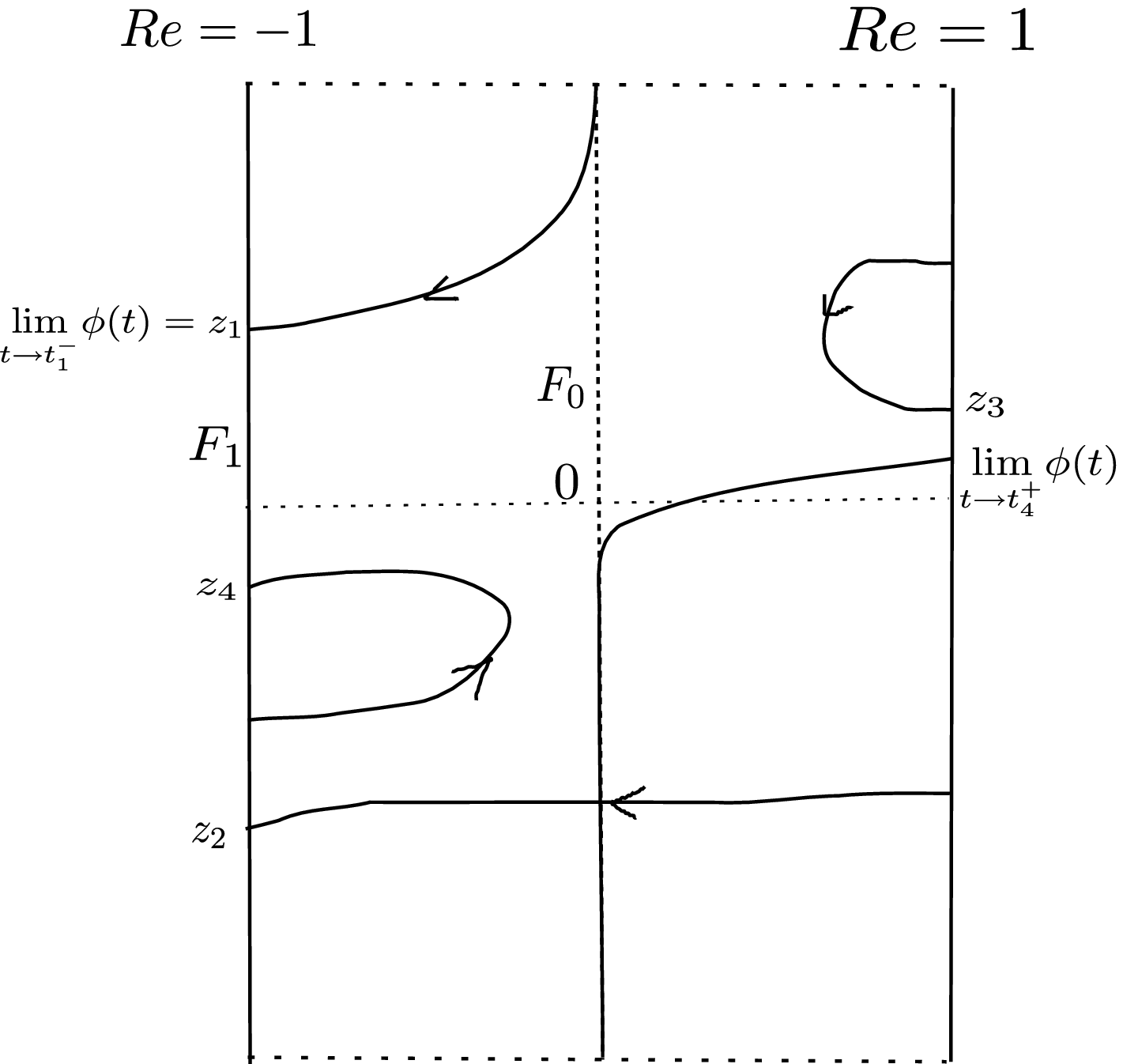}
   \label{fig:subfig1}
   }
 \subfigure[Pushing the curve]{
  \includegraphics[scale=0.4]{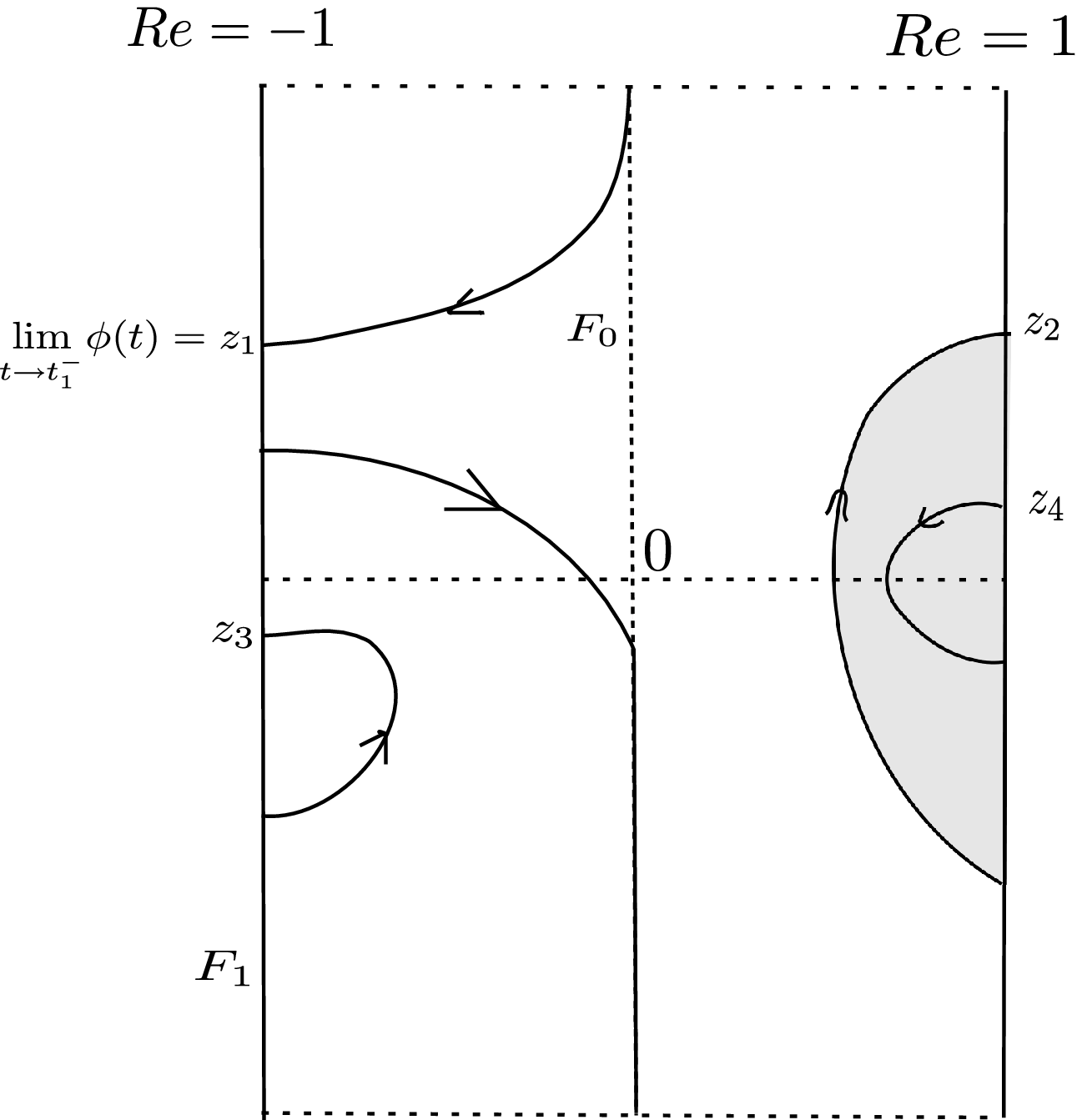}
   \label{fig:subfig2}
   }
\end{figure}
\nono We give a proof of Proposition \ref{prop:Diffrp2}:

\pfone Thinking of $\RR P^2 = S^2 \slash \sim$, where the equivalence relation identifies antipodal points, the action of $SO(3)$ on $S^2$ preserves equivalence classes and thus descends to an action on $\RR P^2$. Therefore $\Diff(\RR P^2)$ contains $SO(3)$ as a subgroup.  We will show that the homogeneous
space  $\Diff(\RR P^2)/SO(3)$ is weakly contractible.  Fixing an $x\in \RR P^2$, first notice that given $f\in \Diff(\RR P^2)$,
 there exists a unique element $\iota_f\in SO(3)$, such that $\iota_f\circ f$ fixes a framing of $x$, or rather, up to homotopy we may assume it
fixes a neighborhood of $x$. For the uniqueness, we observe that the antipodal map on $S^2$ fixes the equivalence class of the north pole, say, but reverses the orientation of a framing. Therefore, as the complement of a ball in $\RR P^2$ is a M\"{o}bius band, we may identify $\Diff(\RR P^2)/SO(3)$ with the compactly supported diffeomorphism group of
the M\"{o}bius band $B$ with the boundary removed, which we denote as $\Diff_c(B)$. This is homotopic to the diffeomorphisms of the closed M\"{o}bius band which fix the boundary.
Below, we identify $B$ with the bundle $\pi:B \to S^1$ with fibers the unit interval.

We fix a fiber $F_0$ over $p_0\in S^1$ and parameterize $F_0$ as
a map
\[
F_0:(-\infty,+\infty)\rightarrow B.
\]
Define
\[
\fF=\{\phi:(-\infty,\infty)\rightarrow B: \phi(t)=F_0(t)\text{ when }|t|>R\text{ for some }R, \phi\text{ is embedded}\}.
\]
Then for $\phi\in\fF$ we have that $\pi \circ \phi$ is a closed loop in $S^1$ with a well defined degree. Given this, we partition $\fF$ as follows:
\[
\fF_i=\{\phi\in\fF: deg(\pi\circ\phi)=i\}.
\]
\lemmaone \label{flem}
$\fF_i$ is connected when $i=0$, and empty except when $i=-1,0,1$.  Curves in $\fF_{1,-1}$ divide $B$ into two components.
\lemmatwo

\pfone Consider the strip $I=\{|Re(z)|\leq1\}$ in $\CC$, then $B$ is
obtained by gluing the two edges of the strip by $z\sim(-z)$ if
$|Re(z)|=1$ (see Figure (a)).  We denote the distinguished fiber in
$B$ obtained from the glued edges by $F_1$. For a curve
$\phi\in\fF$, we may assume that it intersects $F_1$ transversely.
We thus have a finite subset $T$ of $\RR$, such that
$T=\phi^{-1}(F_1)$.  Write $T=\{t_i\}$ where the $t_i$ are in
increasing order.

%Although we identified the two edges $\{Re(z)=1,-1\}$, we still distinguish the intersections of $\phi$ and $F_1$ by a local orientation.
Formally, we now consider $\phi$ as a map $\phi:\RR \setminus T \rightarrow \mathring{I}$, to the interior of $I$ such that
\[
\lim_{t\rightarrow t_i^+}\phi(t)=-\lim_{t\rightarrow t_i^-}\phi(t).
\]
Let $z_i=\phi(t_i)=\lim_{t\rightarrow t_i^-}\phi(t)$.  For $\phi\in\fF_i$ with $i\neq0$, $T$ must be non-empty.

\medskip
{\it Claim.} There exists an isotopy of $\phi$ to a curve $\phi' \in \fF$ with corresponding points $z'_j$ such that either $Re(z'_j)=-1$ or $Re(z'_j)=1$ for all $j$.

\medskip
{\it Proof of Claim.} If $Re(z_i)=-Re(z_{i+1})$ for some $i$ then
the image of $\phi|_{(t_i, t_{i+1})}$ is a curve in $\mathring{I}$
converging at both ends to points on the same edge. It is possible
that the region formed by $\phi|_{(t_i,t_{i+1})}$ and $F_1$ contains
other such loops (see the shaded area of Figure (b)). If
there are no such loops then the region is empty. Hence we can find
a $j$ such that  $Re(z_j)=-Re(z_{j+1})$ and the region formed by
$\phi|_{(t_i,t_{i+1})}$ and $F_1$ is empty.  Now we can perform an
isotopy to remove the intersections $z_j$ and $z_{j+1}$ by pushing
$\phi$ across the region. After such an isotopy the number of
intersection points with $F_1$ will reduce by $2$ and so after a
finite number we must arrive at a curve satisfying our claim.

Given a curve $\phi$ we may now assume that $Re(z_i)=Re(z_j)$ for all $i$, $j$. Without loss of generality suppose that $Re(z_i)=1$ for all $i$. Then if $\phi \in \fF_0$ we see that $\phi$ avoids $F_1$ completely and thus is isotopic to $F_0$. This proves the first statement.

For the second statement, assume that $|T| \ge 2$, that is, there are intersections $z_1$ and $z_2$ with $T_1$. Then we observe that all paths $\phi|_{(t_i,t_{i+1})}$ must lie beneath $\phi|_{(t_1,t_2)}$ for all $i \ge 2$, and thus cannot converge towards $+\infty$ in $I$. This gives a contradiction thus proving the second statement. The final statement is similarly clear.
\pftwo

\corone \label{conncor} The space $\Diff_c(B)$ is connected.
\cortwo

\pfone Indeed, any $f \in \Diff_c(B)$ maps $F_0$ to a path which, like $F_0$ cannot divide $B$. Thus, by Lemma \ref{flem} the image of $F_0$ lies in $\fF_0$ and, moreover, we may assume up to isotopy that $f$ fixes $F_0$, and by a further isotopy a neighborhood of $F_0$ and the complement of a compact set in $B$. But removing a tubular neighborhood of the boundary and $F_0$ from $B$ leaves a set diffeomorphic to a disk, and as diffeomorphisms of the disk are connected, see \cite{smale}, our corollary follows.
\pftwo

Recall that to prove Proposition \ref{prop:Diffrp2} we must show
that  $\Diff_c(B)$ is contractible. Line fields on $B$ are maps from
$B$ to its projectivized unit tangent bundle, where we identify
vectors differing up to sign. We will only consider fields which are
trivial, that is, coincide with fibers of $B$, outside of a compact
set. The bundle is trivialized by the fibers of $B$ and so line
fields are equivalent to maps from $B$ to $S^1$. Let $l_0$ be the
trivial line field tangent to the fibers. The space of sections
$\mathcal{L}_0$ homotopic to $l_0$ is contractible as all such
sections lift to maps to $\RR$ with compact support.

Given an $f \in \Diff_c(B)$ the line field $f_*l_0$ is homotopic to
$l_0$ by Corollary \ref{conncor}. Thus we have a continuous map
$\Diff_c(B) \to {\mathcal L}_0$. There is also an inverse map which
is well defined at least up to homotopy. For this we need the
following claim.\\

\nono\textit{Claim: Line fields in ${\mathcal L}_0$ have no closed
loops.}

\pfone[Proof of the Claim:] We first observe that any closed loops
must project from $B$ to $S^1$ with degree $1$ or $2$ (up to a
sign). This is because line fields lift to line fields on an
annulus, and it is easily seen that only the generating homotopy
class here can admit a closed orbit. Furthermore, if the loop has
degree $2$ then it bounds a compact region $G$ in $B$. Up to an
isotopy we may assume that $G \cap F_0$ is an interval and we have a
return map from $G \cap F_0$ to itself which reverses the two
boundary points. The return map must have a fixed point which
corresponds to a loop of degree $1$.

Next, we observe that the set of line fields ${\mathcal K} \subset {\mathcal L}_0$ which have
a closed loop of degree $1$ is open. This is because the
Poincar\'{e} return map defined on a suitable interval transverse to
a closed loop is orientation reversing, so the fixed point is stable. However the complement of ${\mathcal K}$ in ${\mathcal L}_0$ is also open. Indeed, if a line field has no closed loops then every integral curve must converge to $\partial B$ in both forward and backward time, and by definition coincides with trivial fibers outside of a compact set. This property is preserved under small perturbations of the line field. Hence ${\mathcal K}$ is also closed, and as $l_0$ has no
closed loops and ${\mathcal L}_0$ is contractible, in particular connected, we deduce that our subset must be empty, proving the
claim. \pftwo

Now, starting with a line field in ${\mathcal L}_0$, given the above
claim all integral curves coincide with fibers of $B$ outside of a
compact set (although each end of a curve may correspond to a different fiber). Thus, following these curves we get a orientation
preserving diffeomorphism from $S^1$ (thought of
as the boundary of $B$) to itself with the following properties:
it does not have any fixed points, and squares to identity.\\

\nono\textit{Claim: The space of such diffeomorphisms, denoted as
$D$, is contractible.}

\pfone[Proof of the Claim:] Fix a point $x_0\in S^1$.  Given $f\in
D$, consider $f(x_0)$ which lies in the contractible set
$S^1\backslash \{x_0\}$. Such assignment $D\rightarrow S^1\backslash
\{x_0\}$ is clearly a fibration.

For any choice of $f(x_0)$, the two points $x_0$ and $f(x_0)$
divide $S^1$ into two closed intervals $I_1$ and $I_2$ (including
these two points themselves).  Therefore, $f$ is identified with a
diffeomorphism from $I_1$ to $I_2$ since it is orientation
preserving.  Such diffeomorphisms are further identified with
$\Diff(I_1)$, which is also contractible. (Thinking of the diffeomorphisms as graphs on the interval, it is a convex set.) The claim then follows.
 \pftwo

Notice that deformations in $D$ can be generated by
deformations of the line fields near the boundary of $B$. Therefore
there is a deformation retract from ${\mathcal L}_0$ to line fields
whose integral curves coincide with the same fiber outside of a
compact set. Up to a choice of parameterizing the curves, such line
fields generate elements of $\Diff_c(B)$ by mapping the fibers onto
the corresponding integral curves.
The resulting map, up to homotopy, is an inverse of the natural map $\Diff_c(B) \to {\mathcal L}_0$ described above.
Hence, $\Diff_c(B)$ is homotopic to ${\mathcal L}_0$, which is
contractible, and the proof is complete.
\pftwo

\end{document}